\documentclass[12pt]{amsart}
\usepackage{amsmath,amsfonts,amssymb}

\def\ssum{\mathop{\sum\!\sum}}

\def\le{\leqslant}

\def\ge{\geqslant}

\renewcommand{\mod}{\mathop{\rm{mod}}}

\def\cP{\mathcal{P}}

\numberwithin{equation}{section}

\theoremstyle{definition}

\begin{document}
\title{ A Weaker But Simpler Sieve Inequality}

\author[Friedlander]{John B. Friedlander}

\maketitle

\dedicatory {\sl \quad\quad\quad\quad\quad For the 75th birthday of Roger Heath-Brown}

\medskip

{\bf Abstract:}
  We discuss a cancellation property of sieve weights, one that is applicable to the distribution of primes and almost-primes in very short intervals.\footnote{ 
key words: primes, sieve weights, short intervals}

\medskip

\section {\bf Introduction}

This note is concerned with a particular sieve lemma, returning us a third time to a topic we first studied forty-five years ago! The set-up is as follows:

Let $2\le z \le D$ be reals, let $\cP$ be a set of primes and
\begin{equation*}
  P=P(z)=\prod_{p\in \cP,\, p<z}p.
\end{equation*}

Let $\lambda_d$ denote a sequence of real numbers supported on integers $d,\, 1 \le d \le D, \, d|P(z)$. The sequences in which we are most interested satisfy the useful conditions
\begin{equation}\label{eq:1.1}
 \lambda_1=1, \,\, |\lambda_d| \le \tau_k(d),
 \end{equation} for all $d\le D$ and some positive integer $k$.
Here the $\tau_k$ are the usual divisor functions.

In Friedlander [2]\footnote{The papers [1,2] were submitted to Roger Heath-Brown at the Cambridge Proceedings. Already, at age twenty-eight, he was a member of the editorial board.} and then with Henryk Iwaniec in our book [3], we had given a bound under certain conditions for the sum
\begin{equation}\label{eq:1.2}
W =\sum_{d | P}\, d\biggl(\sum_{\substack{m | P \\ m \equiv 0(\mod d)}}\, 
    \frac{\lambda_m}{m}\biggr)^2 .
\end{equation}
We are in particular thinking of the $\lambda_d$ as being sieve weights so that, in the case of the combinatorial sieves and specifically the beta-sieve weights [3, formula (6.55)], the last condition~\eqref{eq:1.1} holds with $k=1$ and for the Selberg $\Lambda^2$-sieve weights (see (2.6) below) with $k=3$.

It is obvious that, by simply using the trivial bound $|\lambda_d|\le 1$ for all $d$, one obtains the bound $W\ll (\log D)^3$ or, even worse, $(\log D)^5$ if we take $k=3$. However, as these are sieve weights we should expect a good deal of cancellation. But how much? It turns out to be quite a bit. In [3, Lemma 6.18] with H. Iwaniec we gave, for $\lambda_d$ either the upper or lower bound beta-sieve weights, the following bound:

\begin{equation}\label{eq:1.3}
\sum_{d | P}\, d\biggl(\sum_{\substack{m | P \\ m \equiv 0(\mod d)}}\, 
\frac{\lambda_m}{m}\biggr)^2 \ll \prod_{p\mid P}\Bigl(1 - \frac{1}{p}\Bigr)\ .
\end{equation}

This gave a different proof, in a wider range, from that
previously provided in [2]\footnote{which had required a detailed study of the beta-sieve weights, enabled by the foundational work of Iwaniec [4], then only recently published.}. It turned out that this proof in [3] has
a small flaw, pointed out to us by K. Matom\"aki [5]. Our corrected proof with
H. Iwaniec is given, with acknowledgement, in the appendix.

A main motivation for the inequality~\eqref{eq:1.3} is its application to give upper bounds for the number of primes and lower bounds for the number of almost-primes in most short intervals, specifically for almost all intervals in $[1,x]$ having length $g(x)\log x$ where $g(x)\rightarrow \infty $ arbitrarily slowly, intervals as short as could be reasonably expected. See Section 6.10 of [3] for precise statements.
 More recently this has been greatly developed by K. Matom\"aki [6] to show that there are integers with at most two prime factors in almost all such intervals.

\section {\bf A letter from ``Atle''}

The Department at the University of Toronto, Scarborough, recently had occasion to move (with a bit of a push) into a new building. This caused one to go through a lot of old materials and to jettison a major portion of them. A definite keeper was a letter [7] from A. Selberg, dated November 13, 1981, reminding me of something I seem to have long forgotten. That letter\footnote{Brief quotations in this section are from the letter.}, apart from some personal matters, addressed a question that I had posed to him in an earlier letter.

I had asked about obtaining a specific upper bound for the sum
\begin{equation}\label{eq:2.1}
U = \sum_{d | P}\, \varphi(d)\biggl(\sum_{\substack{m | P \\ m \equiv 0(\mod d)}}\, 
    \frac{\lambda_m}{m}\biggr)^2 ,\quad (\varphi =\,\,\, {\rm Euler's\,\, function}),
\end{equation}
a bound that I had conjectured in [1], near the end of Section 3. His response: `` The answer is in the affirmative''.

Provided that the sieve parameter $\beta$ is chosen not too small, one can obtain that conjectured bound:

{\bf Proposition:} With the $\lambda_m$ as in Lemma 6.18, Section 4,

\begin{equation}\label{eq:2.2}
\sum_{d | P}\, \varphi(d)\biggl(\sum_{\substack{m | P \\ m \equiv 0(\mod d)}}\, 
    \frac{\lambda_m}{m}\biggr)^2  \ll \prod_{p\mid P}\Bigl(1 - \frac{1}{p}\Bigr)\ .
\end{equation}

But wait, looking at it now, isn't this the same upper bound as we have in~\eqref{eq:1.3} for an obviously (slightly) larger sum? So already, ever since the paper [2], that has given us a proof of the stronger result. Why bother re-proving this one now from scratch?

We want to belatedly draw attention to the inequality~\eqref{eq:2.2} for two reasons. In the first place, this is all that is actually needed for the application to the short intervals problem. In the estimation of the sum (6.120) of [3] and earlier in [1,2], the count $h\asymp H, (h,d)=1$ had been quickly bounded by $H$ rather than more carefully by $H\varphi (d)/d$.

The second reason is that a comparison of the proof of~\eqref{eq:2.2} with that of~\eqref{eq:1.3} shows that this smaller sum $U$ is really the more natural of the two. The justification for such a claim stems from the following ``easy to show'' (and easily generalized, see~\eqref{eq:2.10}) identity, stated in the letter.

\smallskip
    {\sl   Let $\lambda_d$ be real numbers and
    $\theta_m =\sum_{d|m}\lambda_d$ for $1\le m \le D$. Then}
   \begin{equation}\label{eq:2.3}
      \sum_{d | P}\, \varphi(d)\biggl(\sum_{\substack{m | P \\ m \equiv 0(\mod d)}}\, 
      \frac{\lambda_m}{m}\biggr)^2 = \prod_{p\in \cP }\bigl(1-\frac{1}{p}\bigr)
      \sum_{\delta|P}\frac{\theta_{\delta} ^2}{\varphi (\delta)} .
      \end{equation}

To see this we interchange the order of summation, switch variables to the complementary divisor of $P$ (i.e. $\eta \rightarrow \delta =P/\eta$) and interchange again. We obtain  
\begin{equation*}
  \begin{aligned}
& \sum_{d | P}\, \varphi(d)\biggl(\sum_{\substack{m | P \\ m \equiv 0(\mod d)}}\, 
    \frac{\lambda_m}{m}\biggr)^2 = \ssum_{m_1|P, m_2|P}\frac{\lambda_{m_1}\lambda_{m_2}}{m_1m_2} \sum_{d |(m_1,m_2)}\varphi(d)\\
    & = \ssum_{m_1|P, m_2|P}\frac{\lambda_{m_1}\lambda_{m_2}}{[m_1,m_2]}
    = \frac{1}{P}\ssum_{m_1|P, m_2|P}\lambda_{m_1}\lambda_{m_2}
    \sum_{\substack{\eta|P\\(\eta,[m_1,m_2])=1}}\varphi(\eta)\\
     & = \frac{\varphi(P)}{P}\sum_{\delta|P}\frac{1}{\varphi(\delta)}
    \ssum_{ m_1|\delta, m_2|\delta}\lambda_{m_1}\lambda_{m_2}
    =\prod_{p\in \cP }\bigl(1-\frac{1}{p}\bigr)\sum_{\delta|P}\frac{\theta_{\delta} ^2}{\varphi (\delta)} . 
  \end{aligned}
\end{equation*}

Remark: The neat identity~\eqref{eq:2.3} for the sum $U$ 
should be compared to the rather more complicated treatment
(see the identity~\eqref{eq:3.2}) needed for the sum $W$.
See also the identity for $W$ in (57) of [6]. Perhaps the identity for $W$ closest in spirit to~\eqref{eq:2.3} is that given in [2, pp.382--383] during the course of the original proof of~\eqref{eq:1.3}, namely
\begin{equation*}
      \sum_{d | P}\, d\biggl(\sum_{\substack{m | P \\ m \equiv 0(\mod d)}}\, 
      \frac{\lambda_m}{m}\biggr)^2 = \prod_{p\in \cP }\bigl(1-\frac{1}{p}\bigr)
      \sum_{\alpha|P}\frac{1}{\alpha \varphi(\alpha)}\sum_{\delta|P/\alpha}\frac{\theta_{\delta}(\alpha)^2}{\varphi(\delta)} ,    
\end{equation*}
where $\theta_{\delta}(\alpha)=\sum_{m|\delta}\lambda_{m\alpha}$.

\medskip

Proceeding from~\eqref{eq:2.3} to the proof of the inequality~\eqref{eq:2.2}
``it is easy to show that the quantity 
$\sum_{\delta|P}\theta_{\delta} ^2/\varphi (\delta)$ can be made bounded''.
The function $\theta_{\delta} ^2$ seems rather difficult to approach so we choose
the trivial bound
\begin{equation}\label{eq:2.4}
\theta_{\delta} ^2 \le \tau(\delta)|\theta_{\delta}| ,
\end{equation}
where, since for the beta-sieve weights we have $|\lambda_d| \le 1$, the case $k=1$ of~\eqref{eq:1.1}. Consequently we have, much more quickly and cleanly than for $W$,
\begin{equation}\label{eq:2.5}
\sum_{\delta|P}\frac{\theta_{\delta} ^2}{\varphi (\delta)} \le \sum_{\delta|P}|\theta_{\delta}|\frac{\tau(\delta)}{\delta}\prod_{p|\delta}\Bigl(1+\frac {1}{p-1}\Bigr)\ ,
\end{equation}
which should be compared with~\eqref{eq:3.3}.

From this point onwards we can take the proof of~\eqref{eq:2.2} to run along the very same lines as does that for~\eqref{eq:1.3} in Section 4, beginning there from~\eqref{eq:3.3}. Despite this smaller bound~\eqref{eq:2.5}, as compared to~\eqref{eq:3.3}, we obtain only the same final statement, now with $U$ in place of $W$.

\smallskip

Remark: Use of the inequality~\eqref{eq:2.4} might appear at first
glance too wasteful since the divisor function really could be much larger, were the $\theta_{\delta}$ arbitrary coefficients. However, for us they come from sieve weights and the remaining factor $|\theta_{\delta}|$ on the right side keeps, within the product, the information that $\delta$ behaves as if it were restricted to almost-primes, hence that $\tau (\delta)$ is rather small. Nevertheless, it might be that at this point we have lost an opportunity to choose a somewhat lesser restriction on the size of $\beta$.

  \medskip

 Selberg [7] points out that with suitable choices the bound~\eqref{eq:2.2} also holds for the $\Lambda^2$ upper-bound sieve. Note that the identity~\eqref{eq:2.3} did not depend on the choice of the sieve coefficients $\lambda_d$ and still holds. In the case of the Selberg sieve we no longer have $|\lambda_d|\le 1$, but now the sieve weights take the form (see for example [3, Chapter 7])
\begin{equation}\label{eq:2.6}
\lambda_d =\sum_{[d_1,d_2]=d}\rho_{d_1}\rho_{d_2} ,
\end{equation}
where the sieve constituents $\rho_d$ satisfy $\rho_1=1, |\rho_d|\le 1$ for all $d$ and $\rho_d =0$ for $d>\sqrt D$. Thus, we find that 
$|\lambda_d|\le 3^{\nu(d)}$ with $\nu(d)=\sum_{p|d}1$, hence (1.1) with $k=3$, so that 
\begin{equation}\label{eq:2.7}
  \theta_m\le \sum_{d|m} 3^{\nu(d)}
  =\sum_{j=0}^{\nu(m)}\binom{\nu(m)}{j}3^j
  = (1+3)^{\nu (m)} = \tau^2(m) .
\end{equation}
Hence, noting that here $\theta_m \ge 0$, we can replace~\eqref{eq:2.4} by
\begin{equation}\label{eq:2.8}
\theta_m ^2 \le \tau^2(m)\theta_m
\end{equation}
and the above remark about the divisor function still applies. 

Following Selberg [7] and using for the most part his notation, we define $f$ to be a multiplicative function satisfying $1<f(p)\le p$ on primes and $f'$ the multiplicative function satisfying $f'(p) = f(p) -1$. We specify
the sieve constituents to be of the form
\begin{equation}\label{eq:2.9}
 \rho_d = \mu(d)\frac{f(d)\sum_{\delta\equiv 0 (\mod d)}1/f'(\delta)}{\sum_{\delta}1/f'(\delta)}
\end{equation}
restricted to $\delta\le \sqrt D$, $\delta|P(z)$.
Note that $d\le \sqrt D$, $|\rho_d|\le 1$
always. 

Again we quote: `` ...it is easy to see that in the $\Lambda^2$ case we get
$$
\sum_{\delta}\frac{\tau^2(\delta)\theta_{\delta}}{\varphi(\delta)}= O(1)  
$$
if we take $f(\delta)$ for instance equal to $\frac 32$ for $p<5$ and $=\frac p4$
for $p\ge 5$.''

\medskip

\section {\bf Some Further Identities}

In this brief section we point to some further identities of nature similar to~\eqref{eq:2.3}.
In the first place, staying with Selberg's notation, hence recalling the definition of the multiplicative functions $f, f'$ from the previous paragraph, we have the following result. 

\smallskip

{\bf Proposition:} Let $\lambda_d$ be real numbers and
$\theta_m =\sum_{d|m}\lambda_d$ for $1\le m \le D$.
Let $f, f'$, $1<f(p) $, $f'(p)= f(p)-1$, be multiplicative functions defined on the finite set of primes $\cP$. Then, with $P=\prod_{p\in \cP}p$, 
\begin{equation}\label{eq:2.10}
        \sum_{d | P}\, f'(d)\biggl(\sum_{\substack{m | P \\ m \equiv 0(\mod d)}}\, 
      \frac{\lambda_m}{f(m)}\biggr)^2 = \prod_{p\in \cP }\biggl(1-\frac{1}{f(p)}\biggr)
      \sum_{\delta|P}\frac{\theta_{\delta} ^2}{f'(\delta)} .
\end{equation}

\smallskip

Here,~\eqref{eq:2.3} is simply the special case $f(m)=m$, $f'(m)=\varphi(m)$, needed for sieving an interval. The proof of this more general version is just the same.

We note a significant aspect of the sum on the left-hand side of~\eqref{eq:2.10}. In the previous paragraphs we have taken the coefficients $\lambda_m$ to be sieve weights, whether for the beta or the $\Lambda^2$ sieve. Now, if instead we take them to be the constituents $\rho_m$, as in~\eqref{eq:2.6}, of a $\Lambda^2$ sieve for a sequence ${\mathcal A} = (a_n)$ with density function $1/f$ then the left-hand side of~\eqref{eq:2.10} gives us the coefficient of the main term (the multiple of $\sum_na_n$) in the sieve upper bound. See (7.1), (7.2) in Selberg's paper [8], wherein we find\footnote{More than a century after Brun, there remains a lack of uniformity in sieve notation. Here in [8], though not usually, Selberg uses the notation $u(n)/n$ in place of $1/f(n)$; throughout [3] (see (5.4) there), we use the notation $g(n)$.}
\begin{equation}\label{eq:2.11}
    \sum_{\delta|P}\frac{\lambda_{\delta}}{f(\delta)}=  \sum_{d | P}\, f'(d)\biggl(\sum_{\substack{m | P \\ m \equiv 0(\mod d)}}\, 
      \frac{\rho_m}{f(m)}\biggr)^2 .
      \end{equation}
This formula, expressing that main term coefficient $\sum_{\delta|P}\lambda_{\delta}/f(\delta)$ as a positive definite quadratic form in the constituents, is useful in the optimization among possible choices for the $\Lambda^2$ weights. One can also find, in the same section of [8], further developed extensions of~\eqref{eq:2.11} suitable for optimization of Selberg lower-bound sieve weights.

When Selberg's paper [8] appeared, a decade after the letter [7], it seemed natural to search, without success as it happened, for the identity~\eqref{eq:2.10} or at least of~\eqref{eq:2.3}. On the other hand, one does find a rather different expression for the sieve main term, valid not just for sieves of Selberg type but also in any sieve bound, upper or lower and which is proven by an argument somewhat similar to that for~\eqref{eq:2.10}. Namely, with the $\lambda_{\delta}$ being any sieve weights, or indeed any real numbers, by (3.3) in [8] we have (yes, by coincidence the same label)
\begin{equation}\label{eq:2.12}
    \sum_{\delta|P}\frac{\lambda_{\delta}}{f(\delta)}=  \prod_{p\in \cP }\biggl(1-\frac{1}{f(p)}\biggr)
      \sum_{d|P}\frac{\theta_d}{f'(d)} . 
 \end{equation}
Equating the right-hand sides of~\eqref{eq:2.11} and~\eqref{eq:2.12} we obtain an identity which might appear slightly different from~\eqref{eq:2.10}, specifically
\begin{equation}\label{eq:2.13}
        \sum_{d | P}\, f'(d)\biggl(\sum_{\substack{m | P \\ m \equiv 0(\mod d)}}\, 
      \frac{\rho_m}{f(m)}\biggr)^2 = \prod_{p\in \cP }\biggl(1-\frac{1}{f(p)}\biggr)
      \sum_{\delta|P}\frac{\theta_{\delta}}{f'(\delta)} , 
      \end{equation}
but actually is the very same since $\theta_{\delta}=\sum_{d|\delta}\lambda_d
= (\sum_{m|\delta}\rho_m)^2$. 

\section{ \bf Appendix: Corrected Proof of Opera Lemma 6.18}

Together with Henryk we take this opportunity to show the following material which dates from October 2020 and, apart from the equation labels and slight changes in the acknowledgement, comes from our uncirculated file ``operarev''.

\medskip

{\bf Lemma 6.18:} {\sl Let $2 \le z \le D^{1/(\beta +1)}$, $\cP$ a set of primes, and $P=P(z)$ the product of those primes $p\in \cP$, $p<z$. Let $\lambda_d$ denote either the upper-bound or lower-bound beta-sieve weights of level $D$  with $\beta \ge 8$. Then, we have}
\begin{equation}\label{eq:3.1}
\sum_{d | P}\, d\biggl(\sum_{\substack{m | P \\ m \equiv 0(\mod d)}}\, 
\frac{\lambda_m}{m}\biggr)^2 \ll \prod_{p\mid P}\Bigl(1 - \frac{1}{p}\Bigr)\ ,
\end{equation}
{\sl the implied constant depending only on $\beta$. }

\begin{proof}\footnote{{\bf Acknowledgement.} In October 2020 we were informed
  by Kaisa Matom\"aki [5] of a ``small inaccuracy in the derivation
  of Lemma 6.18'' caused by our neglect of a certain coprimality condition
  and that she  had, in 
the course of her own researches, given a corrected proof. The proof given 
here is our own correction, but reached after receiving that message. 
The corrected 
formulae in the proof are slightly different than the earlier ones (in [3])
but the statement of the lemma remains completely unchanged.} 
Denote by $W$ the sum on the left-hand side of~\eqref{eq:3.1} and by $V$ the product on the right-hand side of~\eqref{eq:3.1}. Writing
$$
\lambda_m=\sum_{ab=m}\mu (a)\theta_b \ ,
$$ 
we compute the inner sum of $W$ as follows (remember that $d|P$): 

$$\sum_{\substack{m | P \\ m \equiv 0(\mod d)}}\, 
\frac{\lambda_m}{m} = \sum_{\substack{ ab| P \\ ab \equiv 0(d)}}\, 
\frac{\mu (a)\theta_b }{ab} =\sum_{b|P}\theta_b \Phi(b)\ ,
$$
where 
\begin{equation*}
\Phi(b)  = \sum_{\substack{a|P, (a,b)=1\\a\equiv 0(d/(b,d))}}\frac{\mu (a)}{ab}
=\frac{\mu(d) f((b,d))}{\varphi(d)\varphi(b)}\, V 
\end{equation*}
and $f(c) =\mu(c)\varphi(c)$. Hence, the left-hand side 
of~\eqref{eq:3.1} is given by  
\begin{equation}\label{eq:3.2}
W= V^2 \sum_{b_1|P}\sum_{b_2|P}\theta_{b_1}\theta_{b_2} \Phi(b_1,b_2)\, ,
\end{equation}
where
\begin{equation*}
\begin{split}
&\Phi(b_1,b_2) =\bigl(\varphi(b_1)\varphi(b_2)\bigr)^{-1}
  \sum_{d|P}f((b_1,d))f((b_2,d))d \varphi(d)^{-2}\\
 & =\bigl(\varphi(b_1)\varphi(b_2)\bigr)^{-1}
  \prod_{p|P}\bigl( 1 +f((b_1,p))f((b_2,p))p(p-1)^{-2}\bigr)\\
& = \bigl(\varphi(b_1)\varphi(b_2)\bigr)^{-1}
  \prod_{\substack{p|P \\ p\nmid b_1b_2}}\bigl( 1 +\frac{p}{(p-1)^2}\bigr)
  \prod_{\substack{p|b_1b_2 \\ p\nmid (b_1b_2)}}\bigl(\frac{-1}{p-1}\bigr)
  \prod_{p|(b_1,b_2)}(p+1) \\
& = \prod_{p|P}\bigl( 1 +\frac{p}{(p-1)^2}\bigr)\prod_{\substack{p|b_1\\ p\nmid b_2}}
  \frac{-1}{p^2-p+1}\prod_{\substack{p|b_2\\ p\nmid b_1}}\frac{-1}{p^2-p+1}
  \prod_{p|(b_1,b_2)} \frac{p+1}{p^2-p+1}\\
& = \mu(b_1)\mu(b_2)\prod_{p|b_1}(p^2-p+1)^{-1}\prod_{p|b_2}(p^2-p+1)^{-1}
\prod_{p|(b_1,b_2)}(p^3+1)\, V_0\ 
\end{split}
\end{equation*}
with 
$$
V_0 = \prod_{p|P} \Bigl(1+ \frac{p}{(p-1)^2}\Bigr)\ .
$$
Note that 
$$
V_0V^2 = \prod_{p|P} \Bigl(1- \frac{1}{p} +\frac{1}{p^2}\Bigr)\asymp V\ . 
$$
By the trivial bound $|\theta(b)|\le \tau(b)$ we obtain  
$$
|W| \le V_0V^2 \sum_{b_1|P}|\theta_{b_1}| \prod_{p|b_1}(p^2-p+1)^{-1}\Psi(b_1)\ ,
$$
where 
\begin{equation*}
\begin{split} 
& \Psi (b_1)  = \sum_{b|P}\tau(b) \prod_{p|b}(p^2-p+1)^{-1}\prod_{p|(b_1,b)}(p^3+1)\\
& = \prod_{p|P} \Bigl(1+\frac{2}{p^2-p+1}\Bigr) \prod_{p|b_1} 
\Bigl(1+2\frac{p^3+1}{p^2-p+1}\Bigr) \Bigl(1+\frac{2}{p^2-p+1}\Bigr)^{-1}\\
& \asymp \prod_{p|b_1}\Bigl(1+2\frac{p^3+1}{p^2-p+1}\Bigr) \asymp
\prod_{p|b_1}\Bigl(1+2p\bigl( 1+ \frac{1}{p}\bigr)\Bigr)\\
& = \tau(b_1)\prod_{p|b_1} \Bigl(p+\frac32\Bigr)\ .
\end{split}
\end{equation*}
Hence, 
\begin{equation}\label{eq:3.3}
W\ll V\sum_{b|P}|\theta_b|\frac{\tau(b)}{b}\prod_{p|b}\Bigl(1+\frac{5}{2p}\Bigr)\ .
\end{equation}
Since $\theta_1=1$ and $\theta_b$ does not change sign for $b \neq 1$ we get 
$$
W\ll V + V\Bigl| \sum_{b|P}\theta_bh(b)\Bigr|\, ,
$$
where $h(b)$ is the multiplicative function with $h(p)=\frac 2p 
\bigl(1+\frac{5}{2p}\bigr)$. Now we return from the $\theta $'s to 
the $\lambda $'s, obtaining 
\begin{equation*}
\begin{split} 
\sum_{b|P}\theta_bh(b) &= \sum_{ac|P}\lambda_ch(ac)
=\sum_{c|P}\lambda_ch(c)\sum_{\substack{ a| P \\ (a,c)=1}}h(a)\\   
& =\sum_{c|P}\lambda_ch(c)\prod_{\substack{p|P\\p\nmid c}}(1+h(p))
=\prod_{p|P}(1+h(p))\sum_{c|P}\lambda_cg(c)\, , 
\end{split}
\end{equation*}
where $g(c)$ is the multiplicative function with $g(p)=h(p)/(1+h(p))$. 
The density function $g$ satisfies the condition [3, (5.38)] of 
dimension $\kappa = 2$. Recall that in our Lemma 6.18 we assume that
$s=(\log D/\log z) \ge \beta + 1$ and $\beta \ge 8$. Therefore the condition
[3, (6.69)] holds with $\kappa =2$ and we can apply the two bounds in
Propostion 6.7 of [3]. Let $\prod$ denote the product of local factors
$(1-g(p))= (1+h(p))^{-1}$ over $p|P$. If $\lambda_c$ are the upper-bound
sieve weights then [3,(6.31)] and [3,(6.73)] imply
$$
\prod \le \sum_c\lambda_c g(c) \le (1+\psi)\prod 
$$
and if $\lambda_c$ are the lower-bound
sieve weights then [3,(6.32)] and [3,(6.74)] imply
$$
(1-\psi)\prod \le \sum_c\lambda_c g(c) \le\prod 
$$
with $\psi$ a positive constant. In either case we obtain 
$$
\sum_{c|P}\lambda_c\, g(c) \ll \prod_{p|P}\,\bigl(1- g(p)\bigr)
=\prod_{p|P}\bigl(1 + h(p)\bigr)^{-1}\ . 
$$
Hence $W\ll V$, completing the proof of Lemma 6.18. 
\end{proof}

\medskip 
\medskip

Department of Mathematics, University of Toronto

Toronto, Ontario M5S 2E4, Canada  \quad (frdlndr@math.utoronto.ca)

\end{document}